          \documentclass{article}
\usepackage[paperwidth=7in, paperheight=10in, margin=.875in]{geometry}

\usepackage{framed}
\usepackage{amsmath,amssymb}
\usepackage{fancybox,graphicx,epsf,color,hyperref}

\renewcommand{\b}[1]{{\bf #1}}

\newcommand{\mmfunz}[5]{
$$ #1 : \left\{ \begin{array}{ccl}
 #2 & \rightarrow & #3 \\
 #4 & \mapsto& #5   \end{array} \right. $$}
\newcommand{\fz}[3]{#1:\, #2 \rightarrow #3}

\renewcommand{\r}[1]{(\ref{#1})}
\newcommand{\bi}{\begin{itemize}}
\newcommand{\ei}{\end{itemize}}
\newcommand{\be}{\begin{enumerate}}
\newcommand{\ee}{\end{enumerate}}

\newcommand{\bd}{\begin{description}}
\newcommand{\ed}{\end{description}}
\renewcommand{\i}{\item}

\newcommand{\bqn}{\begin{eqnarray}}
\newcommand{\eqn}{\end{eqnarray}}
\newcommand{\eqnn}{\nonumber\end{eqnarray}}
\newcommand{\eqnl}[1]{\label{#1}\end{eqnarray}}

\newcommand{\nn}{\nonumber\\}

\newcommand{\ba}[1]{\begin{array}{#1}}
\newcommand{\ea}{\end{array}}

\newcommand{\R}{\mathbb{R}}

\newcommand{\N}{\mathbb{N}}

\newcommand{\bproof}{\begin{proof}}
\newcommand{\eproof}{\end{proof}}
\newtheorem{Theorem}{\bf Theorem}
\newtheorem{lemma}[Theorem]{\bf Lemma}
\newtheorem{corollary}[Theorem]{\bf Corollary}
\newtheorem{definition}[Theorem]{\bf Definition}
\newtheorem{proposition}[Theorem]{\bf Proposition}
\newtheorem{remark}[Theorem]{\bf Remark}
\newenvironment{proof}[1][]{

\noindent{\it Proof#1.}~~}{\hfill$\square$\medskip}

\newcommand{\bt}{\begin{Theorem}}
\newcommand{\et}{\end{Theorem}}
\newcommand{\bl}{\begin{lemma}}
\newcommand{\el}{\end{lemma}}
\newcommand{\bp}{\begin{proposition}}
\newcommand{\ep}{\end{proposition}}
\newcommand{\bc}{\begin{corollary}}
\newcommand{\ec}{\end{corollary}}
\newcommand{\bdeff}{\begin{definition}}
\newcommand{\edeff}{\end{definition}}
\newcommand{\brem}{\begin{remark}}
\newcommand{\erem}{\end{remark}}

\newcommand{\lam}{\lambda}
\newcommand{\al}{\alpha}
\newcommand{\eps}{\varepsilon}

\newcommand{\Id}{\mathrm{Id}}
\newcommand{\Lip}{\mathrm{Lip}}
\newcommand{\Pt}[1]{\left( #1 \right)}
\newcommand{\Pg}[1]{\left\{ #1 \right\}}
\newcommand{\Pq}[1]{\left[ #1 \right] }

\newcommand{\Pabs}[1]{\left| #1 \right|}

\newcommand{\schema}[1]{\b{\sc #1}}

\newcommand{\supp}{\mathrm{supp}}

\newcommand{\weak}{\rightharpoonup}

\renewcommand{\Pr}{\mathcal{P}_c(\R^d)}
\renewcommand{\H}[1][]{{\bf (H#1)}}

\begin{document}
    \title{A traffic flow model with non-smooth metric interaction: well-posedness and micro-macro limit\thanks{This research was partially supported by the European Research Council under the European Union's Seventh
Framework Program (FP/2007-2013) / ERC Grant Agreement n. 257661.}}


        \author{Paola Goatin\thanks{Inria Sophia Antipolis M\'editerran\'ee, France. {\tt paola.goatin@inria.fr}~~\href{http://www-sop.inria.fr/members/Paola.Goatin/}{http://www-sop.inria.fr/members/Paola.Goatin/}}\and{Francesco Rossi\thanks{Aix Marseille Universit\'e, CNRS, ENSAM, Universit\'e de Toulon, LSIS UMR 7296,13397, Marseille, France. {\tt francesco.rossi@lsis.org}~~\href{http://www.lsis.org/rossif/}{http://www.lsis.org/rossif/}}}}





         \pagestyle{myheadings} \markboth{Traffic flow with non-smooth metric interaction}{P. Goatin, F. Rossi} \maketitle

\maketitle

\begin{abstract}

We prove existence and uniqueness of solutions to a transport equation modelling vehicular traffic in which the velocity field depends non-locally 
on the downstream traffic density via a discontinuous anisotropic kernel. The result is obtained recasting the problem in the space of 
probability measures equipped with the $\infty$-Wasserstein distance. We also show convergence of solutions of a finite dimensional system,
which provide a particle method to approximate the solutions to the original problem.

\end{abstract}

{\bf Keywords:} Transport equations, non-local velocity, Wasserstein distance, macroscopic traffic flow models, micro-macro limits.

{\bf AMS: } Primary: 35F25, 35L65; Secondary: 65M12, 90B20.

\section{Introduction}

In this paper, we are interested in studying the macroscopic traffic flow model introduced in \cite{BlandinGoatin} from the point of view of measure transport equations 
in Wasserstein spaces.

Transport equations with non-local velocities have drawn a growing attention in the mathematical community, starting from the Vlasov equation and other models in kinetic theory, see e.g. \cite{bouchut2000kinetic, dobrushin1979vlasov, spohn2012large}. In this context, {\it non-local} means that the velocity at a given point of the space depends not only on the density at that point, but on the density in a whole neighborhood. The first general results of existence and uniqueness for such equations are given by Ambrosio-Ganbo \cite{ambrosio2008hamiltonian}. There, the authors show that Wasserstein distances are key tools to deal with these equations, since vector fields resulting from non-local interactions are Lipschitz continuous with respect to such distances. Several extensions have been proposed since then, including definition of gradient flows \cite{ambrosio2008gradient}, numerical schemes \cite{pedestrian, piccoli2011time}, generalizations to domains with boundary \cite{figalli2010new} and to transport equations with sources \cite{genwass2, genwass}.

Non-local conservation laws have been introduced recently to model a variety of evolution dynamics: besides road traffic models 
\cite{BlandinGoatin, GoatinScialangaRR2015, Herty, Li, SopasakisKatsoulakis}
and crowd motion models \cite{ColomboGaravelloMercier, cristiani2014multiscale, maury2011, maury2008mathematical, toscani2006kinetic}, they are used to describe granular flows~\cite{AmadoriShen},
sedimenatation~\cite{Betancourt}, conveyor belts~\cite{Gottlich} and aggregation phenomena~\cite{JamesVauchelet}.

As far as road traffic is concerned, an advantage in considering non-local mean velocity depending on a weighted mean of the downstream traffic density is represented 
by the consequent finite acceleration, whose unboundedness is one of the drawbacks of Lighthill-Whitham-Richards \cite{LighthillWhitham, Richards} model
and other classical macroscopic models, which allow for speed jumps. This limits their application in connection with consumption and pollution models, 
which heavily rely on acceleration estimation, see for example comments in \cite{blandin_proceedings}.

In our case, the evolution equation for the density $\rho=\rho(t,x)$ of cars on a (infinite) road is given by the following transport partial differential equation (PDE)
\bqn
\begin{cases}
\partial_t\rho+\partial_x\Pt{\rho \, v\Pt{\int_x^{x+\eta} \rho(t,y)w(y-x)\,dy}}=0, & x\in\R,~t>0,\\
\rho(0,x)=\rho_0(x) & x\in\R,
\end{cases}
\eqnl{e-cauchy}
where $v$ is the mean traffic velocity, $\eta>0$ a given parameter, and $w:[0,\eta]\to\R^+$ is a non-increasing Lipschitz weight with $\int_0^\eta w(x)\,dx=1$.
This is intended to model the fact that drivers adapt their velocity depending on the downstream traffic condition, eventually giving more attention to what happens close to 
them than to cars far beyond. In this respect we speak of {\it metric} interaction, opposite to {\it topological} interaction which take into account the ordering of vehicles:
in this case the influence of preceding vehicles takes into account the presence of other vehicles in between \cite[Section 1.1.1.7]{CristianiPiccoliTosin}.

Solutions of \r{e-cauchy} will be defined in the space {$\mathcal{M}(\R)$} of non-negative measures equipped with the Wasserstein distance, and they are to be intended in the weak sense:
\bdeff
A measure $\mu\in C^0([0,T],\mathcal{M}(\R))$ is a weak solution of \r{e-cauchy} if for all $\varphi\in C_c^\infty((0,T)\times\R)$ there holds
\[
\int_0^T\int_\R \left(\partial_t \varphi + \partial_x \phi \, v\Pt{\int_x^{x+\eta} w(y-x)\,d\mu(t,y)}\right) d\mu(t,x) =0.
\]
\edeff

The main result of the paper guarantees existence and uniqueness of weak solutions. 

\bt\label{t-main} Let  $v:[0,1]\to\R^+$ be a Lipschitz and non-increasing function with $v(1)=0$ and $\rho_0\in BV_c(\R,[0,1])$. Then the Cauchy problem \r{e-cauchy}
admits a unique solution in the weak sense, which satisfies
\bqn
0\leq \rho(t,x) \leq \max\Pg{\rho_0}\qquad \mbox{for a.e.~~} x\in \R,t>0.
\eqnn
\et

Above, $BV_c$ denotes the set of functions of bounded variation with compact support.
The proof will be given in Section \ref{s-appr-proof}. Remark that, unlike \cite{BlandinGoatin, GoatinScialangaRR2015}, the proof of uniqueness of solutions does not rely on any entropy condition. 
On the other side, we must restrict the problem to initial data with compact support. 

Observe that the result is based on the study of \r{e-cauchy} in a more general setting, with $v$ not necessarily decreasing. Assuming that $v$ is Lipschitz and $w$ is a (possibly discontinuous) $BV$ interaction kernel of bounded variation, we prove in Proposition \ref{p-gen} existence and uniqueness of the solution of \r{e-cauchy} for small times. It is interesting to observe that the density can blow-up in a finite time $T$, leading to non-existence of the solution for times larger than $T$. Such phenomenon cannot appear for more smooth interaction kernels $w$, see Proposition \ref{p-smooth} and Remark \ref{r-blowup}.

The structure of the article is the following. A short review on Wasserstein distances and general transport equations properties is given in Section~\ref{s-wass}. Transport equations 
with $BV$ interaction kernels are studied in Section~\ref{s-transport}, where we prove existence and uniqueness of solutions locally in time in Proposition~\ref{p-gen}.
 In Section~\ref{s-findim}, we also define a particle approximation for the density, that provides a finite-dimensional numerical scheme for the solution of the transport PDE, and we prove convergence of the micro-macro limit.

\section{Wasserstein distances and transport equations}
\label{s-wass}

\newcommand{\IWD}{$\infty$-Wasserstein distance}
\newcommand{\iw}{W_\infty}
\newcommand{\ess}{\mathrm{ess}}
\renewcommand{\P}{\mathcal{P}}
\newcommand{\Pac}{\mathcal{P}^{ac}}

In this section, we recall the main definitions related to Wasserstein distances, and in particular the definition and properties of the \IWD. We then recall results about transport equations with non-local velocities, in the case of smooth interaction kernels. For more details, see the monographs \cite{old-new,villani} and the articles \cite{pedestrian,genwass,genwass2}. \\

We consider non-negative measures with a given mass $m>0$, which is conserved by the solutions of \r{e-cauchy}. Therefore, without loss of generality, in the following we will deal with probability measures, i.e. $m=1$.

We denote by $\mathcal{P}(\R^d)$ the set of probability measures on $\R^d$ and with $\P_c(\R^d)$ the subset of probability measures with compact support. We also denote with $\Pac(\R^d)$ the subset of probability measures that are absolutely continuous with respect to the Lebesgue measure, and we identify the measure with its density with respect to the Lebesgue measure, e.g. by writing both $\int f(x)\,d\rho(x)$ and $\int f(x)\rho(x)\,dx$ for $\rho\in\Pac(\R^d)$. We use the letters $\mu,\nu,\ldots$ for general measures in $\mathcal{P}(\R^d)$, keeping the notation $\rho,\rho',\ldots$ for measures in $\Pac(\R^d)$. When not specified, we consider $\mathcal{P}(\R^d),\mathcal{P}_c(\R^d),\Pac(\R^d)$ and $\Pac_c(\R^d):=\mathcal{P}_c(\R^d)\cap\Pac(\R^d)$ endowed with the \IWD, defined below. \\

Given a probability measure $\pi$ on $\R^d\times\R^d$, one can interpret $\pi$ as a path to transfer a probability measure $\mu$ on $\R^d$ to another probability measure $\nu$ on $\R^d$ as follows: each infinitesimal mass on a location $x$ is sent to a location $y$ with a probability given by $\pi(x,y)$. Formally, $\mu$ is sent to $\nu$ if the following property holds:
\[ 
\int_{\R^d} d\pi(x,\cdot)=d\mu(x),\qquad \qquad \int_{\R^d} d\pi(\cdot,y)=d\nu(y),
\] 
or, equivantly, for all $f,g\in C^\infty_c(\R^d)$ 
\[
\int_{\R^d\times \R^d} (f(x)+g(y))\,d\pi(x,y) = \int_{\R^d} f(x)\,d\mu(x)+ \int_{\R^d} g(y)\,d\nu(y)
\]
In this case, $\pi$ is called a transference plan from $\mu$ to $\nu$. We denote by $\Pi(\mu,\nu)$ the set of such transference plans. 

Fix now $p\in[1,+\infty)$. One can define a cost for $\pi$ as follows $$J\Pq{\pi}:=\int_{\R^d\times\R^d} |x-y|^p \,d\pi(x,y)$$ and look for a minimizer of $J$ in $\Pi(\mu,\nu)$. Such problem is called the Monge-Kantorovich problem. A minimizer of $J$ in $\Pi(\mu,\nu)$ always exists, see \cite{villani}. A natural space on which $J$ is finite is the space of Borel probability measures with finite $p$-moment, that is
\bqn
\mathcal{P}_p(\R^d):=\Pg{\mu\in\mathcal{P}(\R^d)\ |\ \int |x|^p\, d\mu(x)<\infty}. 
\eqnl{e-moment}
The \b{$p$-Wasserstein distance} is defined on $\mathcal{P}_p(\R^n) \times \mathcal{P}_p(\R^n)$ as
\bqn
W_p(\mu,\nu)=\min_{\pi\in\Pi(\mu,\nu)} J\Pq{\pi}^{1/p}.
\eqnn
It is indeed a distance on $\mathcal{P}_p(\R^d)$, see \cite{villani}. Several topological properties are of interest for the space $\mathcal{P}_p(\R^d)$ endowed with the Wasserstein distance $W_p$, see \cite{old-new,villani}. For future use, here we recall the following.

\bp \label{p-weakp} The Wasserstein distance metrizes weak convergence in $\mathcal{P}_p(\R^d)$, i.e.
$$W_p(\mu_n,\mu)\to 0$$ if and only if $$\mu_n\weak \mu \mbox{~~~~ and ~~~~} \lim_{R\to+\infty}{\lim\sup}_{n\to\infty} \int_{|x|>R} |x|^p\, d\mu_n(x)=0.$$
In particular, $W_p(\mu_n,\mu)\to 0 \mbox{ if and only if } \mu_n\weak \mu$ in $\mathcal{P}(K)$ with $K$ compact.
\ep

\bp \label{p-completep} The space $\mathcal{P}_p(\R^d)$ endowed with the Wasserstein distance $W_p$ is complete.
\ep

\bp \label{p-orderp} The Wasserstein distances are ordered, i.e. $p\leq q$ implies
$$W_p(\mu,\nu)\leq W_q(\mu,\nu).$$
\ep

We finally recall the Kantorovich-Rubinstein duality formula for the 1-Wasserstein distance, see e.g. \cite[Ch. 1]{villani}.
\bp Let $\mu,\nu\in \P_1(\R^d)$. Then it holds
\bqn
W_1(\mu,\nu)=\sup\Pg{\int f(x)\,d(\mu-\nu)(x)\ |\ \Lip(f)\leq 1}.
\eqnl{e-KR}
\ep

\subsection{The $\infty$-Wassersein distance}
In this section, we recall the definition of the \IWD\ and prove some useful properties. We remark that the use of the \IWD\ allows to recover the necessary estimates for the convolution $\int \rho(y) w(y-x)\,dy$ even for BV kernels $w$, see e.g. Proposition \ref{p-C0} below.\\

Given two probability measures $\mu,\nu$ and the space of transference plans $\Pi(\mu,\nu)$ with marginal probabilities $\mu,\nu$, we denote with $C_\infty(\pi)$ the following cost of a transference plan $\pi\in\Pi(\mu,\nu)$:
\bqn
C_\infty(\pi)=\Pg{\pi-\ess\sup(|x-y|)}.
\eqnn
The \IWD\ is then defined as
\bqn
\iw(\mu,\nu)=\inf \Pg{C(\pi)\ |\ \pi\in\Pi(\mu,\nu)}.
\eqnn

We first recall the existence of an optimal transference plan for probability measures with compact support, see \cite[Prop. 2.1]{champion}. Observe that, in analogy with \r{e-moment}, we have $\P_\infty(\R^d)=\P_c(\R^d)$.
\bp \label{p-exist} Let $\mu,\nu\in \mathcal{P}_c(\R^d)$. Then there exists $\pi_*\in\Pi(\mu,\nu)$ realizing $\iw$, i.e.
$$\iw(\mu,\nu)=C_\infty(\pi_*).$$
\ep

It is easy to prove that ordering of the Wasserstein distances is preserved even with the \IWD.
\bp \label{p-order} For any $p\in[1,+\infty)$ it holds
$$W_p(\mu,\nu)\leq \iw(\mu,\nu).$$
\ep
\bproof Let $\pi$ be a transference plan realizing $\iw(\mu,\nu)$. Then it holds
\bqn
W_p^p(\mu,\nu)\leq \int |x-y|^p\,d\pi(x,y)\leq \int C_\infty(\pi)^p \,d\pi(x,y)=\iw^p(\mu,\nu).
\eqnn
\eproof

We now prove lower semicontinuity of $\iw$ with respect to the weak convergence of measures in a compact space.
\bp \label{p-lower} Let $K$ be a compact set in $\R^d$, and $\mu_n$ a sequence in $\mathcal{P}(K)$. If $\mu_n\weak\mu$ and $\nu\in\mathcal{P}_c(\R^d)$, then
$$W_\infty(\nu,\mu)\leq {\lim\inf}_{n\to\infty}\iw(\nu,\mu_n).$$
\ep
\bproof First observe that $\supp(\mu)\subset K$. Since $\supp(\nu)$ is compact, eventually replacing $K$ by $K\cup \supp(\nu)$, one has that all measures have compact support in $K$.

We prove the result by passing to a minimizing subsequence (that we do not relabel) and prove that $W_\infty(\nu,\mu)\leq {\lim}_{n\to\infty}\iw(\nu,\mu_n)$. By Proposition \ref{p-exist}, for each $n$ there exists a transference plan $\pi_n$ realizing $\iw(\nu,\mu_n)$. Since $\pi_n\in\mathcal{P}(K\times K)$ and $K\times K$ is compact, then Prokhorov's Theorem ensures the existence of a subsequence (that we do not relabel) for which it holds $\pi_n\weak \pi_*$ for some $\pi_*\in\mathcal{P}(K\times K)$. It is easy to prove that $\pi_*\in\Pi(\nu,\mu)$. Since $C_\infty$ is lower semicontinuous with respect to the weak topology (see e.g. \cite[Lemma 2.3]{champion}), then $C_\infty(\pi_*)\leq \lim_n C_\infty(\pi_n)$. By recalling that  $\iw(\nu,\mu)\leq C_\infty(\pi_*)$, the result is proved.
\eproof

We now prove results related to the topology of $\mathcal{P}_c(\R^d)$ endowed with the \IWD.
\bp \label{p-complete} The space $\mathcal{P}_c(\R^d)$ is complete with respect to the metric $\iw$.
\ep
\bp \label{p-weak}
Let $\mu_n$ a sequence in $\mathcal{P}_c(\R^d)$. If $\iw(\mu_n,\mu)\to 0$ for some $\mu\in\mathcal{P}_c(\R^d)$, then $\mu_n\weak\mu$.
\ep
\bproof Let $\mu_n$ be a Cauchy sequence in $\mathcal{P}_c(\R^d)$ endowed with the metric $\iw$. For a given $\eps>0$, consider $N$ such that for all $n\geq N$, $k>0$ it holds $\iw(\mu_n,\mu_{n+k})<\eps$. In particular it holds $\iw(\mu_N,\mu_{n})<\eps$, that in turn implies that $\supp(\mu_n)\subset \cup\Pg{ B(x,\eps)\ |\ x\in\supp(\mu_N)}$ for all $n\geq N$. Since such set is bounded, and the supports of $\mu_k$ for $k<N$ are bounded too, then there exists a compact set $K$ containing the support of all $\mu_n$.

Since the $\mu_n$ have uniformly bounded support, then they also have uniformly bounded $p$-th moment for each $p\in[1,\infty)$. Recall that each space $\mathcal{P}_p(\R^d)$ is complete with respect to $W_p$, see Proposition \ref{p-completep}. Observe that it holds $W_p(\mu_n,\mu_m)\leq W_\infty(\mu_n,\mu_m)$, then $\mu_n$ is a Cauchy sequence in $\mathcal{P}_p(\R^d)$, hence there exists $\mu_*$ for which $W_p(\mu_n,\mu_*)\to 0$ 
for all $p\in[1,\infty)$. Then, by Proposition \ref{p-weakp}, we have $\mu_n\weak \mu_*$, proving Proposition \ref{p-weak}.

We are now left to prove that $\iw(\mu_n,\mu_*)\to 0$, that is a direct consequence of lower semicontinuity of $\iw$ with respect to the weak topology of measures, proved in Proposition \ref{p-lower}.
\eproof

\brem It is false that weak convergence of measures implies convergence with respect to the metric $\iw$, even in $\mathcal{P}(K)$ with $K$ compact. For example, consider the following sequence $\mu_n:=\frac{n-1}{n}\delta_0+\frac{1}{n}\delta_1$, that weakly converges to $\mu_*:=\delta_0$, but for which it holds $\iw(\mu_n,\mu_*)=1$ for all $n$. This is in sharp contrast with the $W_p$ metric with $p\in[1,+\infty)$, as recalled in Proposition \ref{p-weakp}.
\erem

\subsection{Transport equations with smooth non-local interactions} \label{s-transportsmooth}
In this section, we study the transport equation with non-local interactions, i.e. the following Cauchy problem:
\bqn
\begin{cases}
\partial_t\mu+\nabla_x\cdot\Pt{V[\mu]\mu}=0, & x\in\R^d, ~t>0,\\
\mu(0,x)=\mu_0(x), &  x\in\R^d,
\end{cases}
\eqnl{e-generale}
where $V$ is a function that associates to each measure $\mu$ a vector field $V[\mu]$. For the simplest case of $V$ actually not depending on $\mu$, the solution of \r{e-generale} is given by the push-forward of the flow of $V$. We recall its definition here.

\bdeff \label{p-pushforward} Given a Borel map $\fz{\gamma}{\R^d}{\R^d}$, the push-forward of a probability measure $\mu\in\P(\R^d)$ is defined by:
\bqn
\gamma\#\mu(A):=\mu(\gamma^{-1}(A))
\eqnn
for every subset $A$ such that $\gamma^{-1}(A)$ is $\mu$-measurable.
\edeff

If $v:\R^d\times[0,T] \to\R^d$ is a vector field uniformly Lipschitz with respect to the space variable and continuous with respect to time, then we denote by $\gamma_t=\phi_t^v$ is the flow generated by $v$, where $\phi_t^v(x_0)$ is the unique solution at time $t$ of 
\[
\begin{cases}
\dot{x}=v(t,x), \\
x(0) = x_0.
\end{cases}
\]

Then, we have the following result, see \cite[Thm 5.34]{villani} and \cite{pedestrian}\footnote{The proof of \r{e-stime4} is given in \cite{pedestrian} for $p<+\infty$, but it can be easily adapted to $p=+\infty$.}.
\bp \label{p-vcost} Let $v:\R^d\times[0,T] \to\R^d$ be a uniformly Lipschitz vector field. Then the equation 
\bqn
\begin{cases}
\partial_t\mu+\nabla_x \cdot\Pt{v\, \mu}=0, & x\in\R^d, ~t>0,\\
\mu(0,x)=\mu_0(x), &  x\in\R^d,\end{cases}
\eqnn
with $\mu_0\in\P_c(\R^d)$ admits a unique solution $\mu\in C^0([0,T],\P_c(\R^d))$. Such solution satisfies $\mu(t)=\phi^v_t\#\mu_0$. 
In particular, if $\mu_0\in\Pac(\R^d)$, then $\mu(t)\in\Pac(\R^d)$ for all times $t>0$.

Moreover, let $v,w$ be two Lipschitz vector fields with Lipschitz constant $L$ and bounded, and $\phi^v_t,\phi^w_t$ the corresponding flows. 
Let $\mu,\nu\in\mathcal{P}_c(\R^d)$ be two probability measures. Then, for $p=[1,+\infty]$ it holds
\bqn
W_p(\phi^v_{t}\#\mu,\phi^w_{t}\#\nu)\leq e^{\frac{p+1}{p} L t} W_p(\mu,\nu)+\frac{e^{L t/p}(e^{Lt}-1)}{L} \sup_{t\in[0,T]}\|v(.,t)-w(.,t)\|_{C^0}.
\eqnl{e-stime4}
\ep

For $V$ actually depending on $\mu$, we have the following theorem, generalizing the results in \cite{ambrosio,pedestrian,genwass,genwass2}.
\bt\label{t-PR} Let $p\in[1,+\infty]$. Let the function
\mmfunz{V\Pq{\mu}}{\mathcal{P}_c(\R)}{(C^{1}\cap L^\infty)(\R^d;\R^d)}{\mu}{V\Pq{\mu}}
satisfy 
\bi
\i $V\Pq{\mu}$ is uniformly Lipschitz and uniformly bounded, i.e. there exist $L$, $M$ not depending on $\mu$, such that for all $\mu\in\Pr, x,y\in\R^d,$
\bqn
\hspace{-5mm}|V\Pq{\mu}(x)-V\Pq{\mu}(y)|\leq L |x-y|\,, \qquad |V\Pq{\mu}(x)|\leq M.
\eqnn

\i $V$ is a Lipschitz function, i.e. there exists $K$ such that 
\bqn\|V\Pq{\mu}-V\Pq{\nu}\|_{\mathrm{C^0}} \leq K W_p(\mu,\nu).
\eqnn
\ei
Then the Cauchy problem \r{e-generale} admits a unique solution $\mu\in C^0([0,T],\mathcal{P}_c(\R^d))$ for all times $T>0$. Moreover, if the initial data $\mu_0$ satisfies $\mu_0\in\Pac(\R^d)$, then $\mu(t)\in\Pac(\R^d)$ for all times $t\in[0,T]$.

Finally, if $\mu,\mu'$ are solutions of \r{e-generale} with initial data $\mu_0,\nu_0$ respectively, there holds
\bqn
W_p(\mu(t),\nu(t))\leq e^{(4L+4K)t} W_p(\mu_0,\nu_0).
\eqnl{e-contdipgen}
\et
\bproof The proof is given for $p<+\infty$ in \cite[Prop. 4 and Thm. 2]{pedestrian} and in \cite{genwass,genwass2}. The key estimate is \r{e-stime4}, that holds also for $p=+\infty$. Then, the original proofs can be easily adapted to $p=+\infty$.
\eproof

We now adapt such result to our setting, in which $d=1$ and $V[\mu](x):=v(\int_\R  w(y-x)\,d\mu(y))$ with $w$ Lipschitz.
\bp\label{p-smooth} Let $v:\R\to\R$ be a Lipschitz and bounded function, and $w:\R\to\R$ be a Lipschitz function with bounded support. If $\mu_0\in \mathcal{P}_c(\R)$, then the Cauchy problem
\bqn
\begin{cases}
\partial_t\mu+\partial_x\Pt{\mu v\Pt{\int_\R w(y-x)\,d\mu(t,y)}}=0, & x\in\R,~t>0,\\
\mu(0,x)=\mu_0(x), & x\in\R,
\end{cases}
\eqnl{e-cauchysmooth}
admits a unique solution $\mu\in C^0([0,T],\mathcal{P}_c(\R))$ for all times $T>0$.

Moreover, if $\mu,\nu$ are solutions of \r{e-cauchysmooth} with initial data $\mu_0,\nu_0$ respectively, there holds
\bqn
W_p(\mu(t),\nu(t))\leq e^{8\Lip(v)\Lip(w)t} W_p(\mu_0,\nu_0),
\eqnl{e-contdip}
for any $p\in[1,+\infty]$.
\ep
\bproof It is sufficient to prove that $V[\mu]$ defined by $V[\mu](x):=v(\int_\R w(z-x)\,d\mu(z))$ satisfies the hypotheses of Theorem \ref{t-PR}. We have
\bqn
\Pabs{V[\mu](x)-V[\mu](y)}&\leq &\Lip(v)\int\Pabs{w(z-x)-w(z-y)}\,d\mu(z)\leq \Lip(v)\Lip(w)|x-y|,\nn
|V[\mu](x)|&\leq& \sup(|v|),\nn
\Pabs{V[\mu](x)-V[\nu](x)}&\leq& \Lip(v)\Pabs{\int w(z-x)\,d(\mu-\nu)(z)}\leq \Lip(v)\Lip(w)W_1(\mu,\nu)\leq\nn
&\leq& \Lip(v)\Lip(w)W_p(\mu,\nu),
\eqnn
where the last estimate is based on the Kantorovich-Rubinstein duality \r{e-KR} and ordering of Wasserstein distances in Propositions \ref{p-orderp} and \ref{p-order}. By identifying $L,K$, we find \r{e-contdip}.
\eproof

\section{Transport equations with BV interactions kernels}
\label{s-transport}
\newcommand{\spazio}{L^\infty\cap\mathcal{P}^{ac}_c(\R)}

In this section, we study transport equations in one space dimension with non-local velocities given by interaction kernels of the form \r{e-cauchy}. 
Unlike Proposition \ref{p-smooth}, here we do not assume that 
the kernel interaction $w$ is Lipschitz continuous, but only $BV$. This prevents to use the results in Section \ref{s-transportsmooth}. 
Moreover, the interaction $\int w(z-x)\,d\mu(z)$ itself is not well-defined for general probability measures, but for measures that are absolutely continuous with respect to the Lebesgue measure, i.e. for $\mu\in\Pac(\R)$ only. 

In particular, we will often deal with measures in the space $\spazio$, endowed with the $\iw$ distance. As stated above, the space $\P_c(\R)$ is complete with respect to the $\iw$ distance, but this does not hold anymore for $\spazio$. Nevertheless, we have convergence if both the support and the $L^\infty$ norm are uniformly bounded, as proved in the following proposition.
\bl \label{l-stack}
Let $K\subset\R$ be a compact set, and $\mu_n\in L^\infty\cap\mathcal{P}^{ac}(K)$ be a sequence of measures with uniformly bounded $L^\infty$ norm and weakly converging to $\mu$. Then $\mu\in L^\infty\cap\mathcal{P}^{ac}(K)$ and $\|\mu\|_{L^\infty}\leq \lim\sup_{n\to\infty}\|\mu_n\|_{L^\infty}$.
\el
\bproof It is clear that $\mu\in\mathcal{P}(K)$. We now prove that $\mu$ is absolutely continuous with respect to the Lebesgue measure $\lam$, by proving that, for each set $A$ that is $\lam$-measurable with $\lam(A)=0$, it holds $\mu(A)=0$. 

We first prove $\mu(A)=0$ in the case of $A$ open. Fix $(a,b)\supset K$ and define $A_n:=\Pg{x\in[a,b]\ |\ d(x,\R\setminus A)\geq \frac1n}$ with $d(x,B)=\inf_{y\in B} |x-y|$. Observe that the $A_n$ are compact, they satisfy $A_n\subset A_{n+1}$ and it holds $\lim_n A_n= A$. Define a sequence of continuous functions $f_n$ with support in $[a-1,b+1]$ satisfying $\chi_{A_n}\leq f_n\leq \chi_A$, increasing in the sense that $f_{n}\leq f_{n+1}$. This implies $f_n\weak \chi_A$.
Hence, by the monotone convergence theorem and weak convergence, we have
\bqn
\mu(A)=\sup_{n\in\N}\int f_n\,d\mu\leq \sup_{n\in\N}{\lim\sup}_{k\to\infty} \int f_n\,d\mu_k\leq \sup_{n\in\N}\int f_n \,dM\lam\leq M\lam(A)=0
\eqnl{e-utile}
with $M\geq \sup_{n\in\N}\|\mu_n\|_{L^\infty}$.

By regularity of finite measures (see \cite{ev-gar}), the result holds for any $\lam$-measurable set $A$, then $\mu\in\Pac(\R)$. The estimate $\|\mu\|_{L^\infty}\leq \lim\sup_{n\to\infty}\|\mu_n\|_{L^\infty}$ is again a consequence of \r{e-utile}.
\eproof\\

The goal of this section is to prove the following result of existence and uniqueness for small times.
\bp \label{p-gen} Let the following hypotheses hold:
\begin{framed}
\begin{description}
\i[\H:] The function $v:[0,1]\to\R$ is Lipschitz and bounded. The interaction kernel $w$ satisfies $w\in BV([\al,\beta],\R^+)$ for fixed $\al,\beta\in\R$, extended with zero in $\R\setminus[\al,\beta]$, and $\int_\al^\beta w(x)\,dx=1$. The initial density $\rho_0$ satisfies $\rho_0\in\spazio$. 
\end{description}
\end{framed}
Then, there exists $T>0$ such that for all $t\in(0,T)$ there exists a unique weak solution $\rho\in C^0([0,t],\spazio)$ of the Cauchy problem
\bqn
\begin{cases}
\partial_t \rho (t,x) + \partial_x(V\Pq{\rho(t)}(x) \rho(t,x))=0, & x\in\R,~t>0,\\
\rho(0,x)=\rho_0(x), & x\in\R.
\end{cases}
\eqnl{e-gen}
where 
$$V\Pq{\rho(t)}(x)=v\Pt{\int \rho(t,y) w(y-x)\,dy}.$$
\ep

To prove this proposition, we first study in the following sections a set of useful technical lemmas related to BV kernels and the corresponding transport equations. The proof of Proposition \ref{p-gen} will be then given in Section \ref{s-proofgen}.

\subsection{Estimates for convolutions with BV kernels}

In this section, we prove some useful estimates for BV functions and functions defined by convolutions with BV kernels.

\bl \label{l-TV} Let $f\in BV([\alpha,\beta],\R)$. Let $[a,b]\subset[\alpha,\beta]$ and $h$ such that $\alpha\leq a-h\leq b+h\leq \beta$. Then it holds
\bqn
\int_a^b \Pabs{f(x)-f(x-h)}\,dx\leq |h|\, TV(f).
\eqnn
\el
Above, $TV(f)$ denotes the (bounded) total variation of $f$.

\bproof  See also \cite[Remark 3.25]{AFP}. Assume that $f$ is non-decreasing and $h>0$. It holds
\bqn
&0&\leq \int_a^b f(x)-f(x-h)\,dx \nn
&&=\int_a^b (f(x)-f(a-h))\,dx-\int_a^b (f(x-h)-f(a-h))\,dx \nn
&&=\int_a^b (f(x)-f(a-h))\,dx-\int_{a-h}^{b-h} (f(x)-f(a-h))\,dx\nn
&& \leq \int_{b-h}^b (f(x)-f(a-h))\,dx\nn
&&\leq \int_{b-h}^b (f(b)-f(a-h))\,dx=h(f(b)-f(a-h)) \nn
&&\leq h\, TV(f).
\eqnn
The proof for $h<0$ is identical.

Let $f=g-l$ with $g,l$ non-decreasing and $TV(f)=TV(g)+TV(l)$. Then it holds
\bqn
&\displaystyle{\int_a^b \Pabs{f(x)-f(x-h)}\,dx}&=\int_a^b \Pabs{\Pt{g(x)-g(x-h)}-\Pt{ l(x)-l(x-h)}}\,dx\nn
&&\leq\int_a^b \Pabs{g(x)-g(x-h)}\,dx+\int_a^b \Pabs{l(x)-l(x-h)}\,dx \nn
&&\leq |h| \Pt{TV(g)+TV(l)} \nn
&&=|h|\,TV(f).
\eqnn
\eproof

\bl\label{l-stimaint} Let $\mu,\mu'\in\spazio$ and $w\in BV([\alpha,\beta],\R)$. Then it holds
\bqn
\Pabs{\int w(x)\,d(\mu(x)-\mu'(x))\,dx}\leq W_\infty(\mu,\mu') \,TV(w)\, \min\Pg{\|\mu\|_{L^\infty},\|\mu'\|_{L^\infty}}
\eqnl{e-eq1}
\el
\bproof Let $\pi\in\Pi(\mu,\mu')$ be a transference plan realizing $h:=\iw(\mu,\mu')$, that exists due to \cite[Prop. 2.1]{champion}. Decompose $w=f-g$ with $f,g$ non-decreasing functions via the Jordan decomposition on the interval $[\alpha,\beta]$. We have
\bqn
&\displaystyle{\int w(x)\,d(\mu(x)-\mu'(x))}&=\int (f(x)-g(x))\,d\pi(x,y)-\int (f(y)-g(y))\,d\pi(x,y) \nn
&&=\int (f(x)-f(y))\,d\pi(x,y)+\int (g(y)-g(x))\,d\pi(x,y) \nn
&&\leq \int (f(x)-f(x-h))\,d\pi(x,y)+\int (g(x+h)-g(x))\,d\pi(x,y) \nn
&&=\int (f(x)-f(x-h))\,d\mu(x)+\int (g(x+h)-g(x))\,d\mu(x) \nn
&&\leq \|\mu\|_{L^\infty} \Pt{\int (f(x)-f(x-h))\,dx+\int (g(x+h)-g(x))\,dx} \nn
&&\leq \|\mu\|_{L^\infty} h \,TV(w),
\eqnn
where we used that points $(x,y)$ in the support of $\pi$ satisfy $|x-y|\leq h$ except for a set of zero measure. By replacing $w$ with $-w$, we have the absolute value on the left hand side of \r{e-eq1}. Since the estimate is symmetric with respect to $\mu,\mu'$, one has the result.
\eproof
\brem The main reason for which Lemma \ref{l-stimaint} holds only for measures $\mu,\mu'$ that are absolutely continuous with respect to the Lebesgue measure is that one needs to give sense to the integral $\int w(x)\,d\mu(x)$ when $w$ is a $BV$ function. The proposition holds for real measures since we need to use the Jordan decomposition of the $BV$ functions.
\erem

\bp \label{p-Lip} Let $w\in BV([\al,\beta],\R)$ and $\rho\in L^\infty(\R)$. Then the function $f(x):=\int \rho(y) w(y-x)\,dy$ is Lipschitz, with Lipschitz constant $L\leq \|\rho\|_{L^\infty} TV(w)$.
\ep
\bproof By using Lemma \ref{l-TV}, we have
\bqn
|f(x_1)-f(x_2)|&\leq& \int \rho(y) |w(y-x_1)-w(y-x_2)|\,dy \nn
&\leq& \|\rho\|_{L^\infty} \int |w(y-x_1)-w(y-x_1+(x_1-x_2))|\,dy \nn
&\leq& \|\rho\|_{L^\infty} |x_1-x_2| \, TV(w).
\eqnn
\eproof
\bp \label{p-C0} Let $w,w'\in BV([\al,\beta],\R)$ and $\rho,\rho'\in L^\infty(\R)$. Then the functions 
$$f(x):=\int \rho(y) w(y-x)\,dy,\qquad f'(x):=\int \rho'(y) w'(y-x)\,dy$$
satisfy
\bqn
\|f-f'\|_{C^0}\leq W_\infty(\rho,\rho') \, TV(w) \min\Pg{\|\rho\|_{L^\infty},\|\rho'\|_{L^\infty}}+\|\rho'\|_{L^\infty} \|w-w'\|_{L^1}.
\eqnn
\ep
\bproof We have
\bqn
\Pabs{f(x)-f'(x)}&\leq& \Pabs{\int \rho(y) w(y-x)- \rho'(y) w(y-x)\,dy}+\Pabs{\int \rho'(y) w(y-x)- \rho'(y) w'(y-x)\,dy}\nn 
&\leq& W_\infty(\rho,\rho') \, TV(w) \min\Pg{\|\rho\|_{L^\infty},\|\rho'\|_{L^\infty}}+\|\rho'\|_{L^\infty} \|w-w'\|_{L^1},
\eqnn
where we used Lemma \ref{l-stimaint} for the first integral and the $L^1$-$L^\infty$ duality for the second integral.
\eproof

\subsection{Solution of transport equations with time-dependent interactions}
In this section, we give estimates for solutions of the transport equations in which the vector field is time-dependent, but not depending on the solution itself.
\bp \label{p-Linf2} Let \H\ hold. Fix $\bar \rho\in C^0([0,T], \spazio)$ and define the time-dependent vector field $V_t$ as follows:
$$V_t(x):=v\Pt{\int \bar\rho(t,y) w(y-x)\,dy}.$$
Define $\rho(t)$ as the unique solution of
\bqn
\partial_t \rho(t,x) + \partial_x(V_t(x) \rho(t,x))=0
\eqnl{e-pde0}
with the given initial data $\rho_0\in\spazio$. Then it holds
\bqn
\iw(\rho(t,\cdot),\rho(t+s,\cdot))\leq s \,\|v\|_{\infty},
\eqnl{e-LipW}
and
\bqn
e^{-Lt}\|\rho_0\|_{L^\infty}\leq \|\rho(t,\cdot)\|_{L^\infty}\leq e^{Lt}\|\rho_0\|_{L^\infty},
\eqnl{e-LipL}
with $L:=Lip(v) \, TV(w) \sup_{s\in[0,t]} \Pg{\|\bar \rho_s\|_{L^\infty}}$.
\ep
\bproof Let $\phi^V_t$ be the flow generated by $V_t$ in the time interval $[t,t+s]$. Then the solution $\rho$ of \r{e-pde0} is unique, and given by the push-forward $\rho(t,\cdot)=\phi_t^V\#\rho_0$, see Proposition \ref{p-vcost}.

We first prove \r{e-LipW}. Since it holds $\rho(t+s,\cdot)=\phi^V_s\#\rho(t,\cdot)$, then the transference plan $\pi(x,y)=(\Id\times \phi^V_s)\#\rho(t,\cdot)$ satisfies $\pi\in\Pi(\rho(t,\cdot),\rho(t+s,\cdot))$. Then it holds
\bqn
\iw(\rho(t,\cdot),\rho(t+s,\cdot))\leq C_\infty(\pi)=\sup_{x\in\supp(\rho(t,\cdot))}\Pg{|x-\phi^V_s(x)|}\leq s\sup (|v|).
\eqnn
We now prove \r{e-LipL}. Observe that the Gronwall lemma gives 
\bqn
e^{-L|t|}|b-a|\leq |\phi^V_t(b)-\phi^V_t(a)|\leq e^{L|t|} |b-a|,
\eqnl{e-gronw}
where $L$ is the Lipschitz constant of $V_s$ in the interval $[0,t]$, that is $L=Lip(v)\, TV(w) \sup_{s\in[0,t]}\Pg{\|\bar \rho_s\|_{L^\infty}}$ 
by Proposition \ref{p-Lip}. 
This implies that for any interval $(a,b)$ there holds
\bqn
\int_{a}^{b} \rho(t,x)\,dx=\int_{a}^{b} \phi^V_t\#\rho_0(x)\,dx=\int^{\phi^V_t(b)}_{\phi^V_t(a)}\rho_0(x)\,dx\leq \|\rho_0\|_{L^\infty} 
|\phi^V_t(b)-\phi^V_t(a)|\leq e^{L|t|} \|\rho_0\|_{L^\infty} |b-a|,
\eqnn
that implies $\|\rho(t)\|_{L^\infty}\leq e^{Lt}\|\rho_0\|_{L^\infty}$. By reversing time, we have the reverse inequality.
\eproof


We now have the following comparison result.
\bp \label{p-key2} 
Let $\bar\rho_t,\bar\rho'_t\in C^0([0,T], \spazio)$ be given and $w,w'$ satisfy \H. Define the  time-dependent vector fields $V_t,V'_t$ as follows:
$$V_t(x):=v\Pt{\int \bar \rho(t,y) w(y-x)\,dy},\qquad V'_t(x):=v\Pt{\int \bar \rho'(t,y) w'(y-x)\,dy}.$$
Define $\rho(t),\rho'(t)$ as the unique solutions of
$$\partial_t \rho(t,x) + \partial_x(V_t(x) \rho(t,x))=0,\qquad\partial_t \rho'(t,x) + \partial_x(V'_t(x) \rho'(t,x))=0,$$
with initial data $\rho_0,\rho'_0$, respectively. Then it holds
\bqn
\iw(\rho(t,\cdot),\rho'(t,\cdot))&\leq& e^{L|t|}\iw(\rho_0,\rho'_0)+\nn
&&+(e^{L|t|}-1)\Pt{\sup_{s\in[0,t]}\iw(\bar\rho(s,\cdot),\bar\rho'(s,\cdot))+\frac{\|w-w'\|_{L^1}}{TV(w')}}\label{e-key2}
\eqn
with $L:=Lip(v)\, TV(w) \sup_{s\in[0,t]}\max\Pg{ \|\bar \rho(s,\cdot)\|_{L^\infty},\|\bar \rho'(s,\cdot)\|_{L^\infty} TV(w')}$.
\ep
\bproof By Proposition \ref{p-Lip}, we have that both $V_t$ and $V'_t$ are Lipschitz vector fields with Lipschitz constant $L$. Denote with $\phi_t,\phi_t'$ the flows of the two vector fields, that are defined for every $t$ with $0\leq t\leq  T$, since the vector fields are globally Lipschitz. Then, we have the following estimate by the Gronwall lemma:
\bqn
\Pabs{\phi_t(x)-\phi_t'(y)}\leq e^{L|t|}|x-y|+\frac{e^{L|t|}-1}{L} \sup_{s\in[0,t]}\|V_s-V'_s\|_{C^0}.
\eqnn
By definition of $V_s$ and $V'_s$, Proposition \ref{p-C0} gives $\|V_s-V'_s\|_{C^0}\leq$
\bqn
 \Lip(v)\Pt{W_\infty (\bar\rho(s,\cdot),\bar\rho'(s,\cdot) \,TV(w)\, \min\Pg{\|\bar\rho(s,\cdot)\|_\infty,\|\bar\rho'(s,\cdot)\|_\infty}+\|\bar\rho'(s,\cdot)\|_{L^\infty} \|w-w'\|_{L^1}}.
\eqnn
By plugging the explicit expression of $L$, and with a careful choice in the $\min$ and $\max$ terms, we have
\bqn
\Pabs{\phi_t(x)-\phi_t'(y)}\leq e^{L|t|}|x-y|+(e^{L|t|}-1) \Pt{\sup_{s\in[0,t]}\iw(\bar\rho(s,\cdot),\bar\rho'(s,\cdot))+\frac{\|w-w'\|_{L^1}}{TV(w')}}.
\eqnn
By reversing time, we also have
\bqn
\Pabs{\phi_t(x)-\phi_t'(y)}\geq e^{-L|t|}|x-y|-(1-e^{-L|t|}) \Pt{\sup_{s\in[0,t]}\iw(\bar\rho(s,\cdot),\bar\rho'(s,\cdot))+\frac{\|w-w'\|_{L^1}}{TV(w')}}.
\eqnn
Take now $\pi_0\in\Pi(\rho_0,\rho'_0)$ and observe that $\pi_t:=(\phi_t\times\phi_t')\#\pi_0$ satisfies $\pi_t\in\Pi(\rho(t,\cdot),\rho'(t,\cdot))$. Consider the set $E_k:=\Pg{(x,y)\in \R^2\ |\ |x-y|>k}$ and observe that it holds
\bqn
\pi_t(E_k)=\pi_0\Pt{\Pg{(\phi_t(x),\phi_t'(y)\ |\ |x-y|>k}}\leq \pi_0\Pt{\Pg{(x,y)\ |\ |x-y|>\tilde k}},
\eqnn
with $\tilde k:=e^{-L|t|}k-(1-e^{-L|t|}) \Pt{\sup_{s\in[0,t]}\iw(\bar\rho(s,\cdot),\bar\rho'(s,\cdot))+\frac{\|w-w'\|_{L^1}}{TV(w')}}$. \\
If $\tilde k>C_\infty(\pi_0)$, one has $\pi_0\Pt{\Pg{(x,y)\ |\ |x-y|>\tilde k}}=0$ by definition of $C_\infty(\pi_0)$. This implies that $\pi_t(E_k)=0$ for all $k$ satisfying
\bqn
k>e^{L|t|} C_\infty(\pi_0)+ (e^{L|t|}-1)\Pt{\sup_{s\in[0,t]}\iw(\bar\rho(s,\cdot),\bar\rho'(s,\cdot))+\frac{\|w-w'\|_{L^1}}{TV(w')}},
\eqnn
hence 
\bqn
C_\infty(\pi_t)\leq e^{L|t|} C_\infty(\pi_0)+ (e^{L|t|}-1)\Pt{\sup_{s\in[0,t]}\iw(\bar\rho(s,\cdot),\bar\rho'(s,\cdot))+\frac{\|w-w'\|_{L^1}}{TV(w')}}.
\eqnl{e-e12}
By recalling that $\iw(\rho(t,\cdot),\rho'(t,\cdot))\leq C_\infty(\pi_t)$ by definition, and passing to the infimum among all transference plans $\pi_0\in\Pi(\rho_0,\rho'_0)$ on the right hand side of \r{e-e12}, we find \r{e-key2}.
\eproof

\subsection{Proof of Proposition \ref{p-gen}} \label{s-proofgen} 

We first prove existence of a solution of \r{e-gen} for small times, via convergence of an explicit Euler scheme. We then prove uniqueness of the solution.

\bp\label{p-ex} Let \H\ hold, and $T>0$ be fixed, with $T<\Pt{e\mathrm{Lip}(v)TV(w)\|\rho_0\|_{L^\infty}}^{-1}$. For each $n\in \N$, consider the following trajectory $\rho^n\in C^0([0,T],\spazio)$:
\bi
\i $\rho^n(0,\cdot):=\rho_0$;
\i $\rho^n(k2^{-n}T+t,\cdot):=\phi_t^{n,k}\#\rho^n(k2^{-n}T,\cdot)$ for $t\in[0,2^{-n}T]$ and $k=0,1,\ldots,2^n-1$, where $\phi_t^{n,k}$ is the flow generated by the vector field 
$$V_{n,k}(x):=v\Pt{\int \rho^n (k2^{-n}T,y) w(y-x)\,dy}.$$
\ei
Then, there exists a subsequence converging to $\rho^*\in C^0([0,T],\spazio)$ that satisfies $\rho^*(0)=\rho_0$ and is a weak solution of 
\bqn
\partial_t \rho(t,x) + \partial_x(V\Pq{\rho(t)}(x) \rho(t,x))=0,
\eqnl{e-pde1}
where 
$$V\Pq{\rho(t)}(x)=v\Pt{\int \rho(t,y) w(y-x)\,dy}.$$
Moreover, $\rho^*$ satisfies $\|\rho^*(t)\|_{L^\infty}\leq e \|\rho_0\|_{L^\infty}$ for all $t\in[0,T]$.
\ep
\bproof Observe that  \r{e-LipW} together with boundedness of $v$ implies that the sequence $\rho^n$ is equibounded and equi-Lipschitz in $C^0([0,T],\mathcal{P}_c(\R))$. Then the Ascoli-Arzel\`{a} theorem implies the existence of a converging subsequence, that we denote again with $\rho^n$. We denote the limit with $\rho^*$, that satisfies $\rho^*(0)=\rho_0$ and $\rho^*\in C^0([0,T],\mathcal{P}_c(\R))$.

Observe that $\mathcal{P}_c(\R)$ is complete with respect to the $\iw$ distance, but this is not the case for $\spazio$. Then, we now prove that $\rho^*(t)\in \spazio$ for $t\in[0,T]$, by proving that $\sup_{t\in[0,T],n\in\N}\|\rho^n(t)\|_{L^\infty}<+\infty$ for $T<\Pt{e\Lip(v)TV(w)\|\rho_0\|_{L^\infty}}^{-1}$ and applying Lemma \ref{l-stack}. Denote  $a^n_k:=\|\rho^n(k2^{-n}T,\cdot)\|_{L^\infty}$ and observe that \r{e-LipL}, together with the definition of $\rho^n((k+1)2^{-n}T,\cdot)$ as function of $\rho^n(k2^{-n}T,\cdot)$, gives the following recurrence rule
$$a^n_{k+1}\leq a^n_{k} e^{2^{-n}T \Lip(v)TV(w) a^n_k},$$ where we have set $a^n_0=\|\rho_0\|_{L^\infty}$. We now prove by induction that $a^n_k\leq e\|\rho_0\|_{L^\infty}$ for all $k\leq 2^n$. We clearly have $a^n_0\leq \|\rho_0\|_{L^\infty}$. Assume that $a^n_j$ satisfies $a^n_j\leq e\|\rho_0\|_{L^\infty}$ for all $j\leq k$. Then we have
\bqn
a^n_{j+1}\leq a^n_{j} e^{2^{-n}T \Lip(v)TV(w) a^n_j}\leq a^n_j e^{2^{-n}}\mbox{~~~~for~~}j=0,\ldots,k,
\eqnn
that in turn implies $a^n_{k+1}\leq a^n_0 e^{k 2^{-n}}\leq e a^n_0$ for $k\leq 2^n$. This implies $\sup_{n\in\N,k=\Pg{0,1,\ldots,2^k}}\|\rho^n(k2^{-n}T,\cdot)\|_{L^\infty}\leq e\|\rho_0\|_{L^\infty}$.

Due to Lipschitzianity of the $L^\infty$ norm given by \r{e-LipL}, we have $$\|\rho^n(t,\cdot)\|_{L^\infty}\leq e^{\Lip(v)e\|\rho_0\|_{L^\infty} TV(w) 2^{-n}T}e\|\rho_0\|_{L^\infty},$$ that in turn implies $\lim\sup_{n\to\infty}\|\rho^n(t,\cdot)\|_{L^\infty}\leq e \|\rho_0\|_{L^\infty}$ for all $t\in[0,T]$. Then, by Lemma \ref{l-stack}, we have $\rho^*\in C^0([0,T],\spazio)$ and $\|\rho^*(t,\cdot)\|_{L^\infty}\leq e \|\rho_0\|_{L^\infty} $ for all $t\in[0,T]$.

We now prove that the limit $\rho^*$ is a weak solution of \r{e-pde1}. Let $\phi\in C^\infty_c((0,T)\times \R,\R)$: we need to prove that 
\bqn
\int_0^T\int_\R\Pt{\rho^*(t,x)\partial_t\phi(t,x)+v\Pt{\int \rho^*(t,y) w(y-x)\,dy}\rho^*(t,x)\partial_x\phi(t,x)}\,dx\,dt=0.
\eqnl{e-weak}
Since by construction of $\rho^n$ the following identity holds
\bqn
\int_0^T\int_\R\Pt{\rho^n(t,x)\partial_t\phi(t,x)+V_{n,k}(x)\rho^n(t,x)\partial_x\phi(t,x)}\,dx\,dt=0,
\eqnn
we prove \r{e-weak} by proving the following three limits:
\bqn
&&\lim_{n\to\infty}\int_0^T\int_\R\Pt{\rho^*(t,x)-\rho^n(t,x)}\partial_t\phi(t,x)\,dx\,dt=0 \label{e-1};\\
&&\lim_{n\to\infty}\int_0^T\int_\R\Pt{v\Pt{\int \rho^*(t,y) w(y-x)\,dy}-V_{n,k}(x)}\rho^*(t,x)\partial_x\phi(t,x)\,dx\,dt=0;\label{e-2}\\
&&\lim_{n\to\infty}\int_0^T\int_\R V_{n,k}(x)\Pt{\rho^*(t,x)-\rho^n(t,x)}\partial_x\phi(t,x)\,dx\,dt=0.
\eqnl{e-3}
Observe that convergence of $\rho^n$ to $\rho^*$ in $C^0([0,T],\spazio)$ implies the existence of a sequence $\eps_n\to 0$ such that $\sup_{t\in[0,T]}\iw(\rho^n(t,\cdot),\rho^*(t,\cdot))<\eps_n$ . Also recall that $\|\rho^*(t,.)\|_{L^\infty}<e \|\rho_0\|_{L^\infty}$ for all $t\in[0,T]$.

To prove \r{e-1}, observe that for each $t\in[0,T]$ the function $\partial_t\phi(t,.)$ is a $BV$ function, hence there exists $C_1:=\sup_{t\in[0,T]}TV(\partial_t\phi(t,.))<+\infty$. Then, by Lemma \ref{l-stimaint}, it holds
\bqn
\int_0^T\Pabs{\int_\R\Pt{\rho^*(t,x)-\rho^n(t,x)}\partial_t\phi(t,x)\,dx}\,dt\leq \int_0^T \eps_n C_1 e \|\rho_0\|_{L^\infty},
\eqnn
that provides \r{e-1}. To prove \r{e-2}, observe that for each $n$ there holds
\bqn
&&\sup_{t\in[0,T],x\in\R}\Pabs{v\Pt{\int \rho^*(t,y) w(y-x)\,dy}-V_{n,k}(x)}\nn
&&\leq\Lip(v) \sup_{k\in\Pg{0,\ldots,2^n -1},t\in[0,2^{-n}T],x\in\R} \Pabs{\int \Pt{\rho^*(k2^{-n}T+t,y)-\rho^n(k2^{-n}T,y)} w(y-x)\,dy}\nn
&&\leq\Lip(v) \sup_{k\in\Pg{0,\ldots,2^n -1},t\in[0,2^{-n}T]} \iw\Pt{\rho^*(k2^{-n}T+t,\cdot),\rho^n(k2^{-n}T,\cdot)} TV(w)\, e\|\rho_0\|_{L^\infty}\nn
&&\leq\Lip(v) \Pt{2^{-n}T \sup(|v|)+\eps_n}TV(w)\, e\|\rho_0\|_{L^\infty},
\eqnl{e-2b}
where we used Lemma \ref{l-stimaint} and the triangular inequality 
\bqn
\iw\Pt{\rho^*(k2^{-n}T+t,\cdot),\rho^n(k2^{-n}T,\cdot)}&\leq&
\iw\Pt{\rho^*(k2^{-n}T+t,\cdot),\rho^*(k2^{-n}T,\cdot)}+\nn
&&\iw\Pt{\rho^n(k2^{-n}T,\cdot),\rho^n(k2^{-n}T,\cdot)}
\eqnn and \r{e-LipW}. 
Then, going back to the left hand side of \r{e-2}, observe that $\|\rho^*(t,\cdot)\|_{L^\infty}\leq e\|\rho_0\|_{L^\infty}$ and integrate on the bounded support of $\partial_x\phi$. By passing to the limit in \r{e-2b}, we have \r{e-2}.

Finally, to prove \r{e-3}, observe that Proposition \ref{p-Lip} provides Lipschitzianity of $V_{n,k}$, with uniform Lipschitz constant $\Lip(|v|)TV(w) e\|\rho_0\|_{L^\infty}$. Also, $\|V_{n,k}\|_{C_0}$ is uniformly bounded by $\|v\|_\infty$. Since one has Lipschitzianity and boundedness of $\partial_x\phi$ too, we have uniform Lipschitzianity of $V_{n,k}\partial_x\phi$. This property, together with boundedness of the support, implies uniform bounded variation, i.e. the existence of $C_2$ such that for all $n,k$ it holds $TV(V_{n,k}\partial_x\phi)<C_2$. Then, Lemma \ref{l-stimaint} provides
\bqn
\int_\R \Pabs{V_{n,k}(x)\Pt{\rho^*(t,x)-\rho^n(t,x)}\partial_x\phi(t,x)}\,dx \leq \iw(\rho^*(t,x),\rho^n(t,x)) C_2\, e\,\|\rho_0\|_{L^\infty}\leq \eps_n C_2\, e\,\|\rho_0\|_{L^\infty}.
\eqnn
By integrating with respect to time and passing to the limit, we have \r{e-3}.
\eproof
\brem\label{r-blowup} It is interesting to observe that one can have blow-up of the density for solutions of \r{e-pde1}. A simple example is given by $w=\chi_{[0,1]}$, $v(x)=x$ and $\rho_0=\chi_{[-1,0]}$. Indeed, observe that the vector field $V\Pq{\rho_0}$ is a non-negative and non-increasing function satisfying $V\Pq{\rho_0}(0)=0$. The evolution $\rho(t)$ has then support contained in $[a(t);0]$ for all times in which it is defined, with $a(0)=-1$ and $a(t)\geq -1$. Observe that $V\Pq{\rho(t)}(a(t))=v(\int_{a(t)}^{a(t)+1}\rho(t)\,dx)=v(1)=1$. Then $a(t)=-1+t$. This implies that $\|\rho(t)\|_{L^\infty}\geq (1-t)^{-1}$, hence $\|\rho(t)\|_{L^\infty}\to \infty$ for $t\to 1$.

This is in sharp contrast with the case of solutions of \r{e-pde1} with Lipschitz kernels $w$. Indeed, in this case hypotheses of Theorem  \ref{t-PR} are satisfied, then one has existence and uniqueness of the solution for all times. Moreover, when the initial data $\mu_0$ is a probability density in $\Pac_c(\R^d)$, then it keeps being in $\Pac_c(\R^d)$.
\erem

We are now ready to prove existence and uniqueness of the solution, i.e. Proposition \ref{p-gen}.

\medskip 

\bproof[ of Proposition \ref{p-gen}] Choose $T<\Pt{e\mathrm{Lip}(v)TV(w)\|\rho_0\|_{L^\infty}}^{-1}$ and apply Proposition \ref{p-ex} to have existence of a solution. Assume now to have two solutions $\bar \rho,\bar \rho'\in C^0([0,T], \spazio)$ and define the time-dependent vector fields $V_t,V'_t$ as follows:
$$V_t(x):=v\Pt{\int \bar \rho(t,y) w(y-x)\,dy},\qquad V'_t(x):=v\Pt{\int \bar \rho'(t,y) w(y-x)\,dy}.$$
Then $\bar \rho,\bar \rho'$ clearly coincide with the solutions $\rho,\rho'$ of 
$$\partial_t \rho(t,x) + \partial_x(V_t(x) \rho(t,x))=0,\qquad\partial_t \rho'(t,x) + \partial_x(V'_t(x) \rho'(t,x))=0,$$
respectively, with initial data $\rho_0$. 

We prove uniqueness by contradiction. Let $t_0:=\inf_{t\in[0,T]}\Pg{\iw(\rho(t,\cdot),\rho'(t,\cdot))>0}$ be the first time in which $\rho,\rho'$ do not coincide. By applying Proposition \ref{p-key2} starting from $t_0$, we have
\bqn
\iw(\rho(t_0+s,\cdot),\rho'(t_0+s,\cdot))\leq (e^{Ls}-1)\sup_{s'\in[0,s]}\iw(\bar\rho({t_0+s'},\cdot),\bar\rho'({t_0+s'},\cdot)),
\eqnl{e-piccolo}
that implies $\sup_{s\in[0,t]}\iw(\rho(t_0+s,\cdot),\rho'(t_0+s,\cdot))\leq (e^{Lt}-1)\sup_{s\in[0,t]}\iw(\bar\rho({t_0+s},\cdot),\bar\rho'({t_0+s},\cdot))$.
By choosing $t<\frac{\ln(2)}{L}$ we have $\iw(\rho(t_0+s,\cdot),\rho'(t_0+s,\cdot))=0$ for all $s\in[0,t]$, that implies 
$$\inf_{t\in[0,T]}\Pg{\iw(\rho(t,\cdot),\rho'(t,\cdot))>0}>t_0,$$ 
which gives a contradiction.
\eproof

Finally, we prove continuous dependence on the interaction kernel and on initial data for the solution of \r{e-pde1}.
\bp\label{p-contdipBV} Let \H\ hold. Let $$T<\Pt{e\mathrm{Lip}(v)TV(w)\max\Pg{\|\rho_0\|_{L^\infty},\|\rho_0\|_{L^\infty}}}^{-1}$$ and $\rho,\rho'$ be the unique solutions of the Cauchy problems
\bqn
\begin{cases}
\partial_t \rho(t,x) + \partial_x(V\Pq{\rho(t)}(x) \rho(t,x))=0,\nn
\rho(0,\cdot)=\rho_0
\end{cases}
&&
\begin{cases}
\partial_t \rho'(t,x) + \partial_x(V'\Pq{\rho'(t)}(x) \rho'(t,x))=0,\nn
\rho'(0,\cdot)=\rho'_0
\end{cases}
\eqnn
where $$V\Pq{\rho(t)}(x)=v\Pt{\int \rho(t,y) w(y-x)\,dy},\qquad V'\Pq{\rho'(t)}(x)=v\Pt{\int \rho'(t,y) w'(y-x)\,dy}.$$
Then it holds
\bqn
\iw(\rho(t,\cdot),\rho'(t,\cdot))\leq e^{4eLt}\iw(\rho_0,\rho'_0)+(e^{4eLt}-1)\frac{\|w-w'\|_{L^1}}{TV(w)},
\eqnl{e-key3}
with $L=\Lip(v)\max\Pg{\|\rho_0\|_{L^\infty} TV(w),\|\rho'_0\|_{L^\infty} TV(w')}$.
\ep
\bproof We apply Proposition \ref{p-key2} with $\rho=\bar\rho$ and $\rho'=\bar\rho'$. Moreover, since $T<\Pt{e\mathrm{Lip}(v)TV(w)\max\Pg{\|\rho_0\|_{L^\infty},\|\rho_0\|_{L^\infty}}}^{-1}$, there holds $\|\rho(t,\cdot)\|_{L^\infty}\leq e\|\rho_0\|_{L^\infty}$, $\|\rho'(t,\cdot)\|_{L^\infty}\leq e\|\rho'_0\|_{L^\infty}$ for all $t\in[0,T]$. Then \r{e-key2} reads as
\bqn
\iw(\rho(t,\cdot),\rho'(t,\cdot))\leq e^{eLt} \iw(\rho_0,\rho'_0)+(e^{eLt}-1)\Pt{\sup_{s\in[0,t]} \iw(\rho(t,\cdot),\rho'(t,\cdot))+\frac{\|w-w'\|_{L^1}}{TV(w)}}
\eqnn
with $L=\Lip(v)\max\Pg{\|\rho_0\|_{L^\infty} TV(w),\|\rho'_0\|_{L^\infty} TV(w')}$. By applying the Gronwall lemma, we have \r{e-key3}.
\eproof

\section{Finite dimensional approximation: the micro-macro limit}
\label{s-findim}

In this section, we describe a finite-dimensional approximation for the solution of \r{e-cauchy}, that may serve as a numerical scheme to compute such solution. 

We first define a Lipschitz approximation $w_\ell$ for the interaction kernel $w$.
\bdeff Let $\eta>0$ be fixed, and $w:[0,\eta]\to\R^+$ be a non-increasing $C^1$ function such that $\int_0^\eta w(x)\,dx=1$. We define $w_\ell:\R\to\R^+$ as follows
\bqn
w_\ell(x):=\begin{cases}
w(0)\frac{\ell+2x}{\ell}&\mbox{~~\hbox{if}~~}x\in\Pq{-\frac{\ell}{2},0}\\
w(x)&\mbox{~~\hbox{if}~~}x\in[0,\eta]\\
w(\eta)\frac{2\eta+\ell-2x}{\ell}&\mbox{~~\hbox{if}~~}x\in\Pq{\eta,\eta+\frac{\ell}{2}}\\
0 &\mbox{~~\hbox{elsewhere.}~~}
\end{cases}
\eqnl{e-wk}
\edeff

We also define a discretization $\Pq{\mu}^n$ of a given probability density $\mu\in \Pac(\R)$. As stated at the beginning, the same idea can be adapted to any density $\mu\in\mathcal{M}^{ac}(\R)$, not necessarily of mass one.
\bdeff \label{d-discr} Let $\mu\in \Pac_c(\R)$ and $n\in\N$ be fixed. Define $\Pg{x_1,\ldots,x_n}$ as follows
\bqn
\begin{cases}
x_1=\sup\Pg{x\in\R\ |\ \int_{-\infty}^x d\mu<\frac{1}{n}},\\
x_{i+1}=\sup\Pg{x\in\R\ |\ \int_{x_i}^x d\mu<\frac{1}{n}},&\mbox{~~~~for $i=1,\ldots,n-1$},
\end{cases}
\eqnl{e-xi}
and 
\bqn
\Pq{\mu}^n:=\frac{1}{n} \sum_{i=1}^n \delta_{x_i}.
\eqnl{e-defrho}
\edeff
Observe that $\Pq{\mu}^n\not\in\Pac(\R)$, hence convolutions with discontinuous kernels are not well defined. We now prove some simple properties of the discretization $\Pq{\mu}^n$.
\bp Let $\mu\in \Pac_c([a,b])$ and $n\in\N$ be fixed, and the $x_i$ and $\Pq{\mu}^n$ defined as in \r{e-defrho}. Then it holds
\bqn
|x_i-x_j|\geq |i-j| \Pt{n\|\mu\|_{L^\infty}}^{-1},
\eqnl{e-ij}
and
\bqn
W_1(\mu,\Pq{\mu}^{n})\leq \frac{b-a}n.
\eqnl{e-W1disc}
\ep
\bproof Assume $i>j$ and observe that it holds
\bqn
\frac{|i-j|}{n}=\int_{x_j}^{x_i}d\mu\leq (x_i-x_j) \|\mu\|_{L^\infty},
\eqnn
that gives \r{e-ij}. To prove \r{e-W1disc}, define $x_0:=a$ and divide the interval $[x_0,x_n]\subset [a,b]$ in the $n$ intervals $[x_0,x_1],[x_1,x_2],\ldots,[x_{n-1},x_n]$. Observing that each of the intervals contains a mass of $\frac{1}{n}$, we have
\bqn
W_1(\mu,\Pq{\mu}^{n})\leq \sum_{i=1}^n W_1\Pt{\mu\chi_{[x_{i-1},x_i]},\frac1n \delta_{x_i}}\leq \sum_{i=1}^n \frac1n \Pt{x_i-x_{i-1}}=\frac{1}{n} \Pt{x_n-x_0}\leq \frac{1}n (b-a).
\eqnn
\eproof

We are now ready to define the approximation $\Pq{\rho}^{n}$ of the solution $\rho$ of \r{e-cauchy}.

\bdeff \label{def-mu} Let $v:[0,1]\to\R^+$ be a Lipschitz and non-increasing function with\footnote{We extend $v$ with the convention $v(x)=0$ for $x>1$.} $v(1)=0$. Let $\eta>0$ be fixed, and $w:[0,\eta]\to\R^+$ be a non-increasing Lipschitz function with $\int_0^\eta w(x)\,dx=1$. Let $\rho_0\in\spazio$ be given.

Fix $n\in\N$ and choose the constant $\ell_n:=(n\|\rho_0\|_{L^\infty})^{-1}$. Define $\Pq{\rho(t)}^{n}$ as 
\bqn
\Pq{\rho(t)}^{n}=\frac{1}{n}\sum_{i=1}^n \delta_{x_i(t)},
\eqnl{empiric}
where $\Pg{x_1(t),\ldots,x_n(t)}$ is the unique solution of the finite-dimension dynamical system
\bqn
\begin{cases}
\dot x_n=v(0),\\
\dot x_i=v\Pt{\frac{1}{n}\sum_{j=1}^n w_{\ell_n}(x_{j}-x_i)}&\qquad\mbox{for $i=1,\ldots,n-1$},\\
x_i(0)=x_{i,0}&\qquad\mbox{for $i=1,\ldots,n$},
\end{cases}
\eqnl{e-findim}
where the $x_{i,0}$ are given by the discretization \r{e-xi} of $\rho_0$, i.e. $\Pq{\rho_0}^n=\frac1n \sum_{i=1}^n \delta_{x_{i,0}}$.
\edeff
It is clear that the solution of \r{e-findim} exists and is unique for all times, due to the Lipschitzianity of $w_{\ell_n}$. Also remark that $\Pq{\rho(t)}^{n}$ is the unique solution of the following transport equation
\bqn
\begin{cases} \label{e-mun}
\partial_t \Pq{\rho(t)}^{n}(x)+\partial_x\Pt{\Pq{\rho(t)}^{n}(x) \, v\Pt{\int w_{\ell_n}(y-x)\,d\Pq{\rho(t)}^{n}(y) }}=0,\\
\Pq{\rho(0)}^{n}=\Pq{\rho_0}^n.
\end{cases}
\eqn
The proof of such identity is direct, by rewriting the interaction kernel $\int w_{\ell_n}(y-x)\,d\mu^n(y)$ in the case of $\mu^n$ given by \r{empiric}.

The main result of this section is the following convergence result, that will be proved in Section \ref{s-appr-proof}.
\bt \label{t-conv} Let  $v:[0,1]\to\R^+$ be a Lipschitz and non-increasing function with $v(1)=0$. Let $\eta>0$ be fixed, and $w:[0,\eta]\to\R^+$ be a non-increasing Lipschitz function with $\int_0^\eta w(x)\,dx=1$. Let $\rho_0\in BV(\R,[0,1])$ with compact support and $\int\rho_0(x)\,dx=1$.

Fix any $T>0$. Then, for each $n\in N$, the approximation $\mu^n$ given by \r{empiric} is defined in $C^0([0,T],\P_c(\R))$ and it satisfies $\Pq{\rho(t)}^{n}\weak \rho(t)$, where $\rho \in C^0([0,T],\P_c(\R))$ is the unique solution of the Cauchy problem
\bqn
\begin{cases}
\partial_t\rho+\partial_x\Pt{\rho\, v\Pt{\int_x^{x+\eta} \rho(t,y)w(y-x)\,dy}}=0, & x\in\R,~t>0,\\
\rho(0,x)=\rho_0(x), & x\in\R .\label{e-problema}
\end{cases}
\eqn
In particular, $\rho$ is defined for all times $t>0$, and it satisfies
\bqn
0\leq \rho(t,x) \leq \max\Pg{\rho_0}\qquad \mbox{for a.e.~~} x\in \R,t>0.
\eqnl{e-maxprinciple}
\et
Before proving Theorem \ref{t-conv}, we address some properties of the finite-dimensional problem \r{e-findim}.

\subsection{Preservation of the minimal and maximal distances}

In this section, we prove that, when the initial minimal distance among agents for the model \r{e-findim} is larger than or equal to $\ell$, then it keeps being larger than or equal to $\ell$ for all times.

\bp\label{p-preservation} Let $\ell\geq\ell_n$ be fixed. Consider a sequence $x_1^0<x_2^0<\ldots<x_n^0$ and denote with $x(t)=(x_1(t),\ldots,x_n(t))$ the unique solution of \r{e-findim}. If $x_{i}^0-x_{i-1}^0\geq \ell$ for all $i=1,\ldots,n$, then it holds $x_{i} (t)-x_{i-1} (t)\geq \ell$ for all times $t>0$.
\ep
\bproof We prove that $x_{i}-x_{i-1}=\ell$ implies $\dot x_{i}-\dot x_{i-1}\geq 0$, that clearly gives the result.

If $i=n$, we have $\dot x_{n}=v(0)=\max_{\sigma\in[0,1]}v(\sigma)\geq \dot x_{n-1}$. For $i<n$, we have
\bqn
\dot x_{i}-\dot x_{i-1}&=&v\Pt{\frac{1}{n}\sum_{j=1}^n w_{\ell_n}(x_{j}-x_{i})}-v\Pt{\frac{1}{n}\sum_{j=1}^n w_{\ell_n}(x_{j}-x_{i-1})}.
\eqnn
Since $v$ is a non-increasing function, we have $\dot x_{i}-\dot x_{i-1}\geq 0$ if and only if
\bqn
\frac{1}{n}\sum_{j=1}^n w_{\ell_n}(x_{j}-x_{i})\leq \frac{1}{n}\sum_{j=1}^n w_{\ell_n}(x_{j}-x_{i-1}) = \frac{1}{n}\sum_{j=2}^{n+1} w_\ell(x_{j-1}-x_{i-1}).
\eqnl{e-qui}
For $j<i$, we have both $x_j-x_i\leq-\ell$ and $x_{j-1}-x_{i-1}\leq-\ell$, thus $w_{\ell_n}(x_j-x_i)=w(x_j-x_i)=0$ 
and $w_{\ell_n}(x_{j-1}-x_{i-1})=w(x_{j-1}-x_{i-1})=0$. 

For $j=i$ we have $w_{\ell_n}(x_{i}-x_i)= w_{\ell_n}(x_{i-1}-x_{i-1}) = w_{\ell_n}(0)= w(0)$. Then, we can restrict ourselves to $j > i$. Condition \r{e-qui} is verified if for each $j=i+1,\ldots,n$ we have 
\bqn
w_{\ell_n}(x_{j}-x_{i})\leq w_{\ell_n}(x_{j-1}-x_{i-1}).
\eqnl{e-condij}  Since $x_{j}-x_{i}>0$ and $x_{j-1}-x_{i-1}>0$ for $j>i$ and $w_{\ell_n}$ is non-increasing in $[0,+\infty]$ as a consequence of the fact that $w$ is non-increasing in $[0,\eta]$, condition \r{e-condij} is equivalent to $x_{j}-x_{i}\geq x_{j-1}-x_{i-1}$, i.e. $x_{j}-x_{j-1}\geq x_{i}-x_{i-1}=\ell$, which is true by hypothesis.
\eproof

The above property is the discrete counterpart of the maximum principle \r{e-maxprinciple}.
The symmetric result is also true, with a similar proof.
\bp\label{p-preservationmax} Let $L\geq\ell_n$ be fixed. Consider a sequence $x_1^0<x_2^0<\ldots<x_n^0$ and denote with $x(t)=(x_1(t),\ldots,x_n(t))$ the unique solution of \r{e-findim}. If $x_{i}^0-x_{i-1}^0\leq L$ for all $i=1,\ldots,n$, then it holds $x_{i} (t)-x_{i-1} (t)\leq L$ for all times $t>0$.
\ep

Observe that, due to the compact support of the solution, the minimum value of the density is always zero, and Proposition \ref{p-preservationmax} 
does not allow to recover a finer minimum principle as in~\cite{BlandinGoatin, GoatinScialangaRR2015}.

\subsection{Convergence to the solution of \r{e-cauchy}}
\label{s-appr-proof}

We now prove Theorem \ref{t-conv}. We combine estimates for solutions of the finite-dimensional problem \r{e-findim} with estimates for solutions of the transport equation with a Lipschitz interaction kernel \r{e-cauchysmooth} and with estimate for solutions of the transport equation with a BV interaction kernel \r{e-cauchy}. 

\bproof[ of Theorem \ref{t-conv}] The idea of the proof is to prove convergence of the approximate solution $\Pq{\rho(t)}^{n}$ to the solution $\rho$ of \r{e-cauchy} by proving intermediate convergence results, divided in four steps. In the first step, we restrict ourselves to a subsequence of $\Pq{\rho(t)}^{n}$ admitting a limit $\rho^*$ in $C^0([0,T],\spazio)$. In the second step, we define a finite-dimensional approximation $\Pq{\tilde\rho(t)}^{n}$ and prove that 
$\Pq{\rho(t)}^{n}$ and $\Pq{\tilde\rho(t)}^{n}$ have the same limit $\rho^*$. In the third step, we define an approximation $\rho^n$ in $\Pac(\R)$ and prove that $\Pq{\tilde\rho(t)}^{n}$ and $\rho^n$ have the same limit $\rho^*$. Finally, in the fourth step we prove that the limit of $\rho^n$ is exactly $\rho$, first for small times and then for any time.

{\bf Step 1.} Fix any $T>0$. We prove that the sequence $\Pq{\rho(t)}^{n}\in C^0([0,T],\P_c(\R))$ admits a subsequence (that we do not relabel) with a limit $\rho^*\in C^0([0,T],\spazio)$, that moreover satisfies $\|\rho^*(t)\|_{L^\infty}\leq \|\rho_0\|_{L^\infty}$. 

Fix $[a,b]$ an interval containing the compact support of $\rho_0$. Then, by construction of $\Pq{\rho(0)}^{n}=\Pq{\rho_0}^n$, we have $\supp(\Pq{\rho(0)}^{n})\subset [a,b]$. Due to boundedness of $v$, we have both $\supp(\Pq{\rho(t)}^{n})\subset [a-T\sup(v),b+T\sup(v)]$ and $\iw(\Pq{\rho(t+s)}^{n},\Pq{\rho(t)}^{n})\leq s \sup(v)$ for all $t\in[0,T]$ and $n\in \N$, i.e. equiboundedness and equi-Lipschitzianity of the sequence $\Pq{\rho(t)}^{n}$ with respect to the $\iw$-distance. Then, eventually passing to a subsequence, there exists a limit $\rho^*\in C^0([0,T],\P_c(\R))$, that has both uniformly bounded support and uniform Lipschitz constant. From now on, the sequence of indexes $n$ is always replaced by the subsequence with limit $\rho^*$.

We now prove that $\|\rho^*(t)\|_{L^\infty}\leq \|\rho_0\|_{L^\infty}$. This is equivalent to prove that for each interval $[\al,\beta]\subset \R$ it holds $\int_\al^\beta \rho^*(t)\,dx\leq \|\rho_0\|_{L^\infty} (\beta-\al)$.  Observe that convergence in $\iw$ implies weak convergence of measures, then $\Pq{\rho(t)}^{n}\weak \rho^*(t)$. For each $\eps>0$, consider a function $\phi_\eps\in C^\infty(\R,\R)$ with support in $[\al-\eps,\beta+\eps]$ such that $\phi_\eps(x)\in[0,1]$ for all $x$ and $\phi_\eps(x)=1$ for $x\in[\al,\beta]$. Observe that, by definition of $\Pq{\rho(t)}^{n}$, for each fixed $t$ it holds 
\bqn
\int  \phi_\eps(x)\, d\Pq{\rho(t)}^{n}(x) =\frac{1}{n}\sum_{x_i(t)\in (\al-\eps,\beta+\eps)}\phi_\eps(x_i(t))\leq \frac{1}{n} \frac{\beta-\al+2\eps}{\ell_n},
\eqnn
where we used the fact that the number of $x_i(t)$ in the interval $(\al-\eps,\beta+\eps)$ is bounded from above, by preservation of the minimal distance $\ell_n=(n\|\rho_0\|_{L^\infty})^{-1}$, see Proposition \ref{p-preservation}. By replacing $\ell_n$ in the last term, by recalling that $\Pq{\rho(t)}^{n}\weak \rho^*(t)$ and observing that $\phi_\eps\geq \chi_{[\alpha,\beta]}$, for each fixed $t$ we have 
\bqn
\int_\al^\beta \,d\rho^*(t,x)\leq \int\phi_\eps(x)\,d\rho^*(t,x)\leq \|\rho_0\|_{L^\infty} (\beta-\al+2\eps).
\eqnn
By passing to the limit for $\eps\to 0$, we have that $\rho^*\in\spazio$ and it satisfies
\bqn
\|\rho^*(t)\|_{L^\infty}\leq \|\rho_0\|_{L^\infty}
\eqnl{e-munLinf}

{\bf Step 2}. We now define $\Pq{\tilde\rho(t)}^{n}$. For each $n$, consider the approximated kernel $w_{m_n}$ with  $m_n:=\ln(n)^{-1}$. Then, define $\Pq{\tilde\rho(t)}^{n}$ similarly to $\Pq{\rho(t)}^{n}$, as follows:
\bqn
\Pq{\tilde\rho(t)}^{n}=\frac{1}{n}\sum_{i=1}^n \delta_{y_i(t)},
\eqnn
where $\Pg{y_1(t),\ldots,y_n(t)}$ is the unique solution of the finite-dimension dynamical system
\bqn
\begin{cases}
\dot y_i=v\Pt{\frac{1}{n}\sum_{j=1}^n w_{m_n}(y_{j}-y_i)}& \mbox{for $i=1,\ldots,n$}\\
y_i(0)=x_{i,0} & \mbox{for $i=1,\ldots,n$},
\end{cases}
\eqnl{e-findim2}
and the $x_{i,0}$ are given by the discretization \r{e-xi} of $\rho_0$, i.e. $\Pq{\rho_0}^n=\frac1n \sum_{i=1}^n \delta_{x_{i,0}}$. Remark that, similarly to \r{e-mun}, $\nu^n$ is the unique solution of the following transport equation
\bqn
\begin{cases}
\partial_t \Pq{\rho(t)}^{n} (x)+\partial_x\Pt{\Pq{\rho(t)}^{n}(x)\, v\Pt{\int w_{m_n}(x-y)\,d\Pq{\rho(t)}^{n}(y))}}=0,\\\label{e-nun}
\Pq{\rho(0)}^{n}=\Pq{\rho_0}^n.
\end{cases}\eqn

We now prove that, for each $T>0$ it holds $\lim_{n\to\infty} \sup_{t\in[0,T]} \iw(\Pq{\rho(t)}^{n},\Pq{\tilde\rho(t)}^{n})=0$. This implies that the limit of both $\Pq{\rho(t)}^{n}$ and $\Pq{\tilde\rho(t)}^{n}$ is $\rho^*$. For each $n$, we define $\eps_i^n(t):=\sup_{s\in[0,t]}|x_i(s)-y_i(s)|$ and $\eps^n(t):=\max_{i=1,\ldots,n}\eps_i^n(t)$. Observe that we have $\iw(\Pq{\rho(t)}^{n},\Pq{\tilde\rho(t)}^{n})\leq \eps^n(t)$, by choosing the transference plan sending $x_i(t)$ to $y_i(t)$. Also observe that $\eps_i^n(0)=\eps^n(0)=0$, since the initial data coincide.

Observe that the initial data satisfies $x^0_{i+1}-x^0_i\geq \ell_n$, then Proposition \ref{p-preservation} for \r{e-findim} with kernel $w_{\ell_n}$ gives $x_{i+1}(t)-x_i(t)\geq \ell_n$. The property does not hold for $y_i$ for big $n$, since $m_n>\ell_n$.

We now estimate the evolution of $\eps_i^n$. By the definition, we have $\eps_n^n=0$. By boundedness of $v$, we have $\eps_i^n(t+s)-\eps_i^n(t)\leq 2\sup(v) s$, hence the $\eps_i^n$ are Lipschitz function. By definition of the dynamics \r{e-findim}, we have
\bqn
\dot\eps_i^n\leq \Lip(v)\frac{1}{n}\sum_{j=1}^n \Pabs{w_{\ell_n}(x_j(t)-x_i(t))-w_{m_n}(y_j(t)-y_i(t))}.
\eqnl{e-epsi}
It is clear that $y_j(t)-y_i(t)=\eps_j^n(t)-\eps_i^n(t)+x_j(t)-x_i(t)$, hence $y_j(t)-y_i(t)\in[x_j(t)-x_i(t)-2\eps^n(t),x_j(t)-x_i(t)+2\eps^n(t)]$. For any fixed index $i$ and for each time $t$, we divide the indexes $j$ in the following five sets:
\bi
\i $j$ is such that $x_j(t)-x_i(t)\in\left]-\infty,-\frac{m_n}{2}-2\eps^n(t)\right]$. In this case, both $x_j(t)-x_i(t)$ and $y_j(t)-y_i(t)$ are in the interval $\left]-\infty,-\frac{m_n}{2}\right]$, for which it holds $w_{\ell_n}(x_j(t)-x_i(t))=w_{m_n}(y_j(t)-y_i(t))=0$.
\i $j$ is such that $x_j(t)-x_i(t)\in\Pt{-\frac{m_n}{2}-2\eps^n(t),2\eps^n(t)}$. In this case, we simply estimate 
$$\Pabs{w_{\ell_n}(x_j(t)-x_i(t))-w_{m_n}(y_j(t)-y_i(t))}\leq w(0).$$ 
Observe that $\Pabs{x_j(t)-x_i(t)}\geq |j-i| \ell_n$, hence the number of indexes $j$ in this set is smaller or equal than $\frac{\frac{m_n}{2}+4\eps^n(t)}{\ell_n}$.
\i $j$ is such that $x_j(t)-x_i(t)\in [2\eps^n(t),\eta-2\eps^n(t)]$: in this case both $x_j(t)-x_i(t)$ and $y_j(t)-y_i(t)$ belong to the interval $[0,\eta]$, on which $w_{\ell_n}$ and $w_{m_n}$ coincide with $w$, that is Lipschitz. Hence $\Pabs{w_{\ell_n}(x_j(t)-x_i(t))-w_{m_n}(y_j(t)-y_i(t))}\leq 2\Lip(w)\eps^n (t)$. In this case the number of indexes $j$ in this set is strictly smaller than $\frac{\eta}{\ell_n}$.
\i $j$ is such that $x_j(t)-x_i(t)\in (\eta-2\eps^n(t),\eta+\frac{m_n}{2}+2\eps^n(t))$. Similarly to the second case, we estimate $\Pabs{w_{\ell_n}(x_j(t)-x_i(t))-w_{m_n}(y_j(t)-y_i(t))}\leq w(0)$ and observe that the number of indexes $j$ in this set is smaller or equal than $\frac{\frac{m_n}{2}+4\eps^n(t)}{\ell_n}$.
\i $j$ is such that $x_j(t)-x_i(t)\in \left[\eta+\frac{m_n}{2}+2\eps^n(t),+\infty\right[$. Similarly to the first case, we have $w_{\ell_n}(x_j(t)-x_i(t))=w_{m_n}(y_j(t)-y_i(t))=0$.
\ei
By using the previous decomposition of the indexes in \r{e-epsi}, we have
\bqn
\dot \eps^n_i(t)&\leq& \Lip(v)\frac{1}{n}\Pt{2w(0) \frac{\frac{m_n}{2}+4\eps^n(t)}{\ell_n}+ 
2\Lip(w)\eps^n (t)\frac{\eta}{\ell_n}}.
\eqnn
Taking the supremum over $\eps^n_i$, and recalling that $\ell_n=(n\|\rho_0\|_{L^\infty})^{-1}$, we have
\bqn
\dot\eps^n(t)&\leq& C_0 m_n+C_1\eps^n(t)
\eqnn
with $C_0:=\Lip(v) w(0) \|\rho_0\|_{L^\infty} $, $C_1:=\Lip(v)\|\rho_0\|_{L^\infty}(8 w(0)+2\Lip(w)\eta)$. Since $\eps^n(0)=0$ and the constants $C_0,C_1$ do not depend on $t$, the Gronwall estimate gives
\bqn
\eps^n(t)\leq \frac{C_0}{C_1} m_n e^{C_1 t}=\frac{C_0}{C_1\ln(n)}e^{C_1 t}.
\eqnn
This implies $W_\infty(\Pq{\rho(t)}^{n},\Pq{\tilde\rho(t)}^{n})\leq \frac{C_0}{C_1\ln(n)}e^{C_1 T}$ for all times $t\in[0,T]$. Then, the limit of $\Pq{\tilde\rho(t)}^{n}$ is $\rho^*$.

{\bf Step 3.} Choose again $[a,b]$ any interval containing the support of the initial datum $\rho_0$, and fix any $T>0$. We define $\rho^n\in C^0([0,T],\spazio)$ as the solution of \r{e-problema} with initial data $\rho_0$ when replacing $w$ with $w_{m_n}$. Since $w_{m_n}$ is Lipschitz, we have existence and uniqueness of $\rho^n$ for any $t\in[0,T]$ thanks to Proposition \ref{p-smooth}. 

Fix now $T_1<(32\Lip(v) w(0))^{-1}$. We now compare $\Pq{\tilde\rho(t)}^{n}$ and $\rho^n(t)$ on the interval $[0,T_1]$, by observing that they are the solution of the same equation 
\bqn
\partial_t \mu+\partial_x\Pt{\mu \, v\Pt{\int w_{m_n}(y-x)\,d\mu(y)}}=0,
\eqnl{e-step21}
with different initial data, that are $\Pq{\tilde\rho(0)}^{n}=\Pq{\rho_0}^n$ and $\rho^n(0,\cdot)=\rho_0$. Since \r{e-step21} satisfies the hypotheses of Proposition \ref{p-smooth}, we have
\bqn
W_1(\Pq{\tilde\rho(t)}^{n},\rho^n(t))&\leq& e^{8\Lip(v)\Lip(w_{m_n}) T_1} W_1(\Pq{\rho_0}^n,\rho_0)\leq  e^{16\Lip(v) w(0) \ln(n)T_1} W_1(\Pq{\rho_0}^n,\rho_0)\leq\nn
&\leq& \sqrt{n} \frac{b-a}{n}=\frac{b-a}{\sqrt{n}}
\eqnn
for a sufficiently big $n$, for which it holds $\Lip(w_{m_n})\leq 2\frac{w(0)}{m_n}\leq 2 w(0)\ln(n)$. We also used \r{e-W1disc} to estimate the initial distance. Then, $\Pq{\tilde\rho(t)}^{n}$ and $\rho^n$ have the same limit in the interval $[0,T_1]$.

To prove the result for the initial interval $[0,T]$, subdivide it in intervals $[kT_1,(k+1)T_1]$ with $k=0,1,\ldots,\frac{T}{T_1}-1$. Then, by induction one can prove $W_1(\Pq{\tilde\rho(t)}^{n},\rho^n(t))\leq \frac{32\Lip(v)\,w(0)\, T(b-a)}{\sqrt{n}}$. As a consequence, the limit of $\rho^n$ is $\rho^*$ on the whole interval $[0,T]$.

{\bf Step 4.} We now fix $T<\Pt{e\mathrm{Lip}(v)TV(w)\|\rho_0\|_{L^\infty}}^{-1}$ and study $W_\infty(\rho^{n}(t),\rho(t))$. In this case, we observe that both $\rho^n$ and $\rho$ are solutions of an equation of the form \r{e-cauchy} with different interaction kernels $w_{m_n}$ and $w$ respectively, and the same initial data $\rho_0$. Observe that, by the particular structure of $w$ and of its approximation $w_{m_n}$, it holds $TV(w)=TV(w_{m_n})=2 w(0)$ and $\|w_{m_n}-w\|_{L^1}\leq \frac{m_n w(0)}2$. Then, by Proposition \ref{p-contdipBV}, we have
\bqn
\iw(\rho^n(t),\rho(t))\leq (e^{4eLt}-1)\frac{\|w_{m_n}-w\|_{L^1}}{TV(w)}\leq (e^{4eLT}-1)\frac{m_n}{4}\leq \frac{e^{4}-1}{4\ln(n)}, 
\eqnn
with $L=2w(0) \Lip(v)\|\rho_0\|_{L^\infty} $. This implies that $\rho=\rho^*$, i.e. convergence of $\Pq{\rho(t)}^{n}$ to $\rho$ in the interval $[0,T]$.

Going back to Step 1, we recall that we chose a converging subsequence of $\Pq{\rho(t)}^{n}$, and we proved that the limit of such sub-sequence is $\rho$, which is  the unique solution of \r{e-cauchy}. Then we have that the whole sequence $\Pq{\rho(t)}^{n}$ converges to $\rho$, in the time interval $[0,T]$ with $T<\Pt{e\mathrm{Lip}(v)TV(w)\|\rho_0\|_{L^\infty}}^{-1}$.

We now prove that $\rho^*=\rho$ in the time interval $[0,T]$ for any $T>0$, and in particular that a solution for \r{e-cauchy} exists (and is unique) for all times. For $T>0$ fixed, divide the interval $[0,T]$ in intervals $[kT',(k+1)T']$ with $T'<\Pt{e\mathrm{Lip}(v)TV(w)\|\rho_0\|_{L^\infty}}^{-1}$ and $k=0,1,\ldots, \frac{T}{T'}-1$. Then, the previous results show that $\rho^*=\rho$ in $[0,T']$. This also implies 
\bqn
\|\rho(T')\|_{L^\infty}= \|\rho^*(T')\|_{L^\infty}\leq \|\rho_0\|_{L^\infty},
\eqnl{e-rhoLinf}
where we used estimate \r{e-munLinf}. Observe that $T'$ satisfies hypotheses of Proposition \ref{p-gen} taking $\rho(T')$ as initial datum: this gives existence and uniqueness of the solution of \r{e-cauchy} on the time interval $[T',2T']$, that coincides with $\rho^*$ on the same interval, and for which it holds $\|\rho(2T')\|_{L^\infty}\leq \|\rho_0\|_{L^\infty}$. By induction, we have $\rho^*=\rho$ on the whole interval $[0,T]$, that also gives $\|\rho(t)\|_{L^\infty}\leq \|\rho_0\|_{L^\infty}$ for all $t\in[0,T]$, i.e. condition $\rho^*(t,x)\leq \max(\rho_0)$ a.e. for $t\in[0,T],x\in\R$.
\eproof




  \bibliographystyle{SIAM}

  \bibliography{nonlocal}

\end{document}